\documentclass[letterpaper, 10 pt, conference]{ieeeconf}  
                                                          
\IEEEoverridecommandlockouts                              
                                                          
\usepackage{amsmath} 
\usepackage{amssymb}  

\usepackage{graphicx}          
\usepackage{mathtools}
\usepackage{accents}
\usepackage{xcolor}

\newcommand{\idf}{\prescript{I}{}{f}}
\newcommand{\iSigma}{\prescript{I}{}{\Sigma}}
\newcommand{\imu}{\prescript{I}{}{\mu}}
\newcommand{\ouf}{\prescript{o}{}{f}}
\newcommand{\oc}[1]{\accentset{\circ}{#1}}

\newcommand{\R}{\mathbb{R}}
\newcommand{\N}{\mathcal{N}}

\newtheorem{rem}{Remark}
\newtheorem{Theorem}{Theorem}
\newtheorem{Definition}{Definition}

\newtheorem{Proposition}{Proposition}

\title{\LARGE \bf Sensitivity and safety of  fully probabilistic control}

\author{Bernat Guillen Pegueroles$^{1,\ast}$ and Giovanni Russo$^{2,\ast}$
\thanks{$^{1}$ Princeton University, Princeton
        {\tt\small bernatp@princeton.edu}}%
\thanks{$^{2}$ University College Dublin
        {\tt\small giovanni.russo1@ucd.ie}}%
        \thanks{$^{\ast}$ Work done in part while the authors were with IBM Research}%
}

\begin{document}

\maketitle
\thispagestyle{empty}
\pagestyle{empty}
\begin{abstract}
In this paper we present a sensitivity analysis for the so-called fully probabilistic control scheme. This scheme attempts to control a system modeled via a probability density function (pdf) and does so by computing a probabilistic control policy that is optimal in the Kullback-Leibler sense. Situations where a system of interest is modeled via a pdf naturally arise in the context of neural networks, reinforcement learning and data-driven iterative control. After presenting the sensitivity analysis, we focus on characterizing the convergence region of the closed loop system and introduce a safety analysis for the scheme. The results are illustrated via simulations.\newline
{\em This is the preliminary version of the paper entitled "On robust stability of  fully probabilistic control with respect to data-driven model uncertainties" that will be presented at the 2019 European Control Conference.}
\end{abstract}


\section{Introduction}

Over the past few years, much research attention has been devoted to the design of model-free and iterative control algorithms that are able to learn how to control a given system of interest by learning from past iterations, see e.g. \cite{8039204} and references therein for a recent application of iterative learning to Model Predictive Control. At the same time, for many cyber-physical systems, driven by the recent {\em explosion} in the amount of available data (via e.g. the Internet of Things) and by the dramatic improvements in computational and communication infrastructures, {\em Deep Learning} techniques, and in particular neural networks, have been increasingly used to model and classify systems \cite{Goodfellow-et-al-2016}. In practice, the output of a deep neural network is often a probability density function (pdf), describing the state of a given system. 

In this context, designing control algorithms for {\em uncertain} systems and subject to certain safety and convergence constraints is rapidly becoming a major research topic. This is motivated by a number of applications, including reinforement learning \cite{Sutton:1998:IRL:551283} and model-free control schemes subject to state constraints \cite{NIPS2017_6692} and with convergence guarantees, see e.g. \cite{Recht2018ATO} for an overview.

The design of control strategies in the absence of reliable models and in the presence of strong uncertainty has long been the subject of stochastic control, see e.g. \cite{nla.cat-vn2527712}. Stochastic control is deeply related to decision science and, in particular, to Bayesian dynamic decision making, see e.g. \cite{Ber_85}, where the control action is computed by minimizing the expected value of a loss function embedding the control goal. The fully probabilistic control algorithm, which was originally introduced in the seminal work \cite{KARNY19961719}, belongs to the family of stochastic control algorithms, with the main difference that it selects randomized control laws that make the entire joint
distribution of closed-loop variables as close as possible (in the sense of the Kulback-Leibler divergence) to their desired distribution. See also \cite{Karny:2005:OBD:1051758,doi:10.1002/acs.742,karny_fully_2006} for recent developments on this topic, together with \cite{HERZALLAH2015199,doi:10.1002/asjc.1717} for the development of an adaptive critic fully probabilistic control, which is closely related to reinforcement learning schemes based on actor critic \cite{6313077}.

In this paper we present a sensitivity and safety analysis for the fully probabilistic control, when the pdfs of the closed loop system are all multivariate Gaussians. In particular, we analyze how error models can propagate through the system and provide a method, based on numerical continuation, to compute the corresponding region of convergence. The region of convergence is a region of the parameter space in which the stability of the closed loop system is guaranteed via the fully probabilistic control. Then, we move onto defining a notion of safety for the system. This notion is related to the maximum error between the target state and the current state of the system. Given this definition and a safety requirement on the closed loop system, we show how it is possible to numerically compute a safety region. This is a region of the parameter space where the norm of the state never violates the safety requirement. Finally, we show how embedding learning mechanisms in the closed loop system can indeed extend both the convergence and the safety regions. The results are illustrated via an example. The complete proofs of the technical results, extended to a non-Gaussian setting, will be presented elsewhere. Also, we do not discuss here the computational aspects of the approach presented in this paper. This aspect will be discussed elsewhere.

\section{Problem set-up}

\subsection*{Notation and initial remarks}

The notation used in this paper is closely related to the one of \cite{karny_fully_2006}.  Let  $S^{\ast}$ be a given finite set, we denote by $\oc{S}$ its cardinality. The value of a given quantity, say $q$, at time $t$ is denoted by $q_t$  and the set $ t^\ast := \{1, \dots, \oc{t}\}$  is a given time horizon. Recall that, given the probability space $(\Omega,\Sigma,\mathbb{P})$ (where $\Omega$ is the sample space, $\Sigma$ is the collection of all the events, and $\mathbb{P}$ is the probability measure) a random variable is a measurable function $X:\Omega\rightarrow \R$ and we denote by $\mathbb{E}[X]$ the expected value of $X$. Also, we denote by $f(\cdot | \cdot)$ a given conditional probability density function (pdf). Then, in the context of this paper, a {\em system} is specified in the probabilistic sense, i.e. the time evolution of the system is specified via the pdf $f(x_t | u_t, x_{t-1})$, where: (i) $x_t$ is the $n$-dimensional observed state of the system at time $t$; (ii) $u_t$ is the $m$-dimensional control input at time $t$. As usual, we denote the multivariate normal distribution of the random vector $v = [v_1,\ldots,v_n]^T$ by $v\sim\N(\mu,\Sigma)$, where $\mu$ is the mean  vector and $\Sigma$ is the covariance matrix. In what follows, $x(t)$ is the sequence of observed states up to time $t$, i.e. $x(t) := (x_1,\ldots,x_t)$ and $u(t)$ is the sequence of observed inputs up to time $t$, i.e. $u(t) = (u_1, \dots, u_t)$. Also, we define the {\em system dataset} (up to time $t$) as $d(t) = (x_0, u_1, x_1, \dots, u_t, x_t)$. Finally, we recall that, given two pdfs, say $f_1$ and $f_2$, over the same set, say $S^\ast$, the Kullback-Leibler (KL) divergence $D(f_1 \| f_2)$ \cite{KL_51} is defined as
$$
D(f_1 \| f_2) := \int f_1(S^\ast) \ln \left(\frac{f_1(S^\ast)}{f_2(S^\ast)}\right)dS^\ast.
$$
We also recall that, see e.g. \cite{Goodfellow-et-al-2016} and references therein, the KL divergence measures the proximity of the pair of pdfs $f_1$ and $f_2$.
\begin{rem}
Applications where the evolution of a system of interest is described via a pdf naturally arise in the context of artificial intelligence and machine learning. For example, certain Bayesian reinforcement learning schemes (such as Thompson Sampling) keep a joint posterior distribution for the belief of the model. Then, the policy is derived to select optimal actions with respect to this posterior \cite{Ghavamzadeh:2015:BRL:2858995.2858996}. Another example is provided by closed loop systems where the underlying process is classified via a neural network. See e.g. \cite{8357977,8317888,7795769} for an application in the context of intelligent transportation. \newline
\end{rem}

\subsection*{Mathematical background}

Assume a system model (specified in the probabilistic sense) is given, together with the {\em ideal} probability distribution
\begin{equation}\label{eq:jointid}
\idf(d(\oc{t}), x(\oc{t})|x_0) = \prod_{t\in t^{\ast}} \idf(x_t | u_t, x_{t-1})\idf(u_t | d(t-1)),
\end{equation}
denoted in what follows by $\idf(\oc{t})$ for simplicity. The ideal distribution is the pdf corresponding to the desired behavior of the system. As shown in \cite{karny_fully_2006}, the pdf can be constructed so as to embed both a given set of control goals and constraints. Indeed, in (\ref{eq:jointid}) the pdf $\idf(x_t | u_t, x_{t-1})$ specifies the desired state evolution and $\idf(u_t | d(t-1))$ specifies the constraints on the evolution of the control law over time. In \cite{karny_fully_2006} an algorithm has been proposed with the goal of producing a {\em control} distribution, $\ouf(u_t | d(t-1))$, $t\in t^{\ast}$, such that the joint probability distribution
\begin{equation}\label{eq:joint}
f(d(\oc{t}) | x_0) = \prod_{t\in t^{\ast}} f(x_t | u_t, x_{t-1})\ouf(u_t | d(t-1)),
\end{equation}
(denoted for simplicity by $f(\oc{t})$) minimizes the Kullback-Leibler divergence $D(f(\oc{t}) \| \idf(\oc{t}))$. In what follows we simply say that $\ouf(u_t|d(t-1))$ is an admissible control strategy (or policy) for the system. The main theoretical result for the design of the control strategy, which has been formalized in \cite{karny_fully_2006} (see Proposition $2$), can be stated as follows:
\begin{Theorem}\label{thm:control}
The optimal admissible control strategy $\ouf(u_t|d(t-1))$ minimizing $D(f(\oc{t}) \| \idf(\oc{t}))$ is the randomized control strategy:
\begin{equation}\label{eq:out}
  \ouf(u_t | d(t-1)) = \frac{\idf(u_t | d(t-1))e^{-\omega(u_t, d(t-1))}}{\gamma(d(t-1))},
\end{equation}    
where, starting from $\gamma(d(\oc{t})) = 1$, $\gamma(\cdot)$ and $\omega(\cdot,\cdot)$ are computed via the following system of backward recursive equations:
 \begin{subequations}\label{eq:control}
\begin{equation}\label{eq:gamma}
\gamma(d(t-1))=\int\idf(u_t|d(t-1))e^{-\omega(u_t, d(t-1))} du_t,
\end{equation}
\begin{equation}\label{eq:omega}
\begin{split}
&\omega(u_t, d(t-1)) =\\
& \int f(x_t | u_t, x_{t-1}) \log\left(\frac{f(x_t|u_t, x_{t-1})}{\gamma(d(t))\idf(x_t|u_t,x_{t-1})}\right)dx_t.
\end{split}
\end{equation}
\end{subequations}
\end{Theorem}

\subsection*{Problem statement}
Theorem \ref{thm:control} is obtained by assuming that the pdf of the model is  known to the control. Unfortunately, this assumption cannot be satisfied in certain applications. For example, in the context of Deep Learning, neural networks can only provide an estimate of the underlying pdf dominating the process that is being observed. Motivated by this, in this paper we present a sensitivity analysis for the above randomized control strategy to errors in the system model, i.e. errors in $f(x_t | u_t, x_{t-1})$. In particular, by focusing on the case where the pdf of the model and the ideal pdf are both Gaussian, we first  perform an analysis of model error propagation and then we draw some considerations on the {\em safety regions} of the closed loop system (this concept is introduced in Section \ref{sec:safety}) and on how iterative learning techniques can be integrated in the control strategy so as to extend the safety region. 

\section{Results}\label{sec:linear}

In this section we present the main results of the paper. We also illustrate the results with a number of simulations.

\subsection{Randomized policy: the Gaussian case}

Assume that: (i) the pdf of the system is Gaussian, in particular $f(x_t| u_t,x_{t-1}) = \N\left(A_tx_{t-1}+B_tu_t,\Sigma_{\Xi}\right)$; (ii) the ideal distributions are also Gaussian distributions, i.e. $\idf(x_t | u_t, x_{t-1}) = \mathcal{N}(\imu_{t,x}, \iSigma_{t,x})$ and $\idf(u_t|d(t-1)) = \mathcal{N}(\imu_{t,u}, \iSigma_{t, u})$. With the following proposition an explicit formula is given for the control strategy (see also \cite{KARNY19961719}).\newline

\begin{Proposition}
	\label{lem:full-linear}
	The optimal admissible control strategy $\ouf(u_t|d(t-1))$ minimizing $D(f(\oc{t}) \| \idf(\oc{t}))$ is the randomized control strategy given by the  system of recursive equations
\begin{subequations}\label{eq:lin}
\begin{equation}\label{eq:linout}
\ouf(u_t|d(t-1)) = \mathcal{N}(\mu_{t, u}, \Sigma_{t, u}),
\end{equation}    
\begin{equation}\label{eq:lingamma}
\begin{split}
\gamma(d(t-1)) &= \\
&C e^{-\frac{1}{2}\left((L_{t-1} x_{t-1} + M_{t-1})^T\Sigma_{t-1, \gamma}^{-1}(L_{t-1} x_{t-1} + M_{t-1})\right)},
\end{split}
\end{equation}
\end{subequations}
with $\gamma(d(\oc{t})) = 1$, $C$ being a constant not depending on the data and where the control parameters are given by the following set of backward recursion equations:

\begin{subequations}\label{eq:linear_rec}
\begin{equation}\label{eq:sig-u-t} 
\Sigma_{t,u}^{-1} = \iSigma_{t,u}^{-1} + B_t^T \Sigma_{t,\omega}^{-1} B_t,
\end{equation}    
\begin{equation}\label{eq:mu-u-t} 
\begin{split}
&\mu_{t,u} = \\
&\Sigma_{t,u}(\iSigma_{t,u}^{-1}\imu_{t,u} + B_t^T\Sigma_{t,\omega}^{-1}(\mu_{t,\omega}-A_tx_{t-1})),
\end{split}
\end{equation}
\begin{equation}\label{eq:l-t} 
L_{\oc{t}} = 0, \ \ L_{t-1} = A_{t}, \ \  t \leq \oc{t},
\end{equation}
\begin{equation}\label{eq:sig-g-t} 
\Sigma_{t-1,\gamma}^{-1} = (B_t \iSigma_{t,u} B_t^T + \Sigma_{t,\omega})^{-1}, \ \ \Sigma_{\oc{t}, \gamma}^{-1} = 0,
\end{equation}
\begin{equation}\label{eq:m-t} 
\begin{split}
& M_{t-1} = \\
&\Sigma_{t-1, \gamma}(\Sigma_{t, \omega}^{-1}B_t\Sigma_{t,u}\iSigma_{t,u}^{-1}\imu_{t,u} - \mu_{t, \omega}), \ \  M_{\oc{t}} = 0,
\end{split}
\end{equation}
\begin{equation}\label{eq:sig-om-t} 
\Sigma_{t,\omega}^{-1} = \iSigma_{t, x}^{-1} + L_t^T \Sigma_{t, \gamma}^{-1}L_t,
\end{equation}
\begin{equation}\label{eq:mu-om-t} 
\mu_{t, \omega} = \Sigma_{t, \omega}(\iSigma_{t,x}^{-1}\imu_{t,x} - L_t^T\Sigma_{t,\gamma}^{-1}M_t).
\end{equation}\newline
\end{subequations}
\end{Proposition}

{\em Sketch of the proof.} The proof is obtained by induction. First, the hypotheses imply that  $\ouf(u_t|d(t-1))$ is also Gaussian, i.e.
$$
\ouf(u_t|d(t-1)) = \mathcal{N}(\mu_{t, u}, \Sigma_{t, u}).
$$
This gives (\ref{eq:linout}) and we will sketch how $\mu_{t, u}, \Sigma_{t, u}$ can be devised from the recursive equations (\ref{eq:linear_rec}). Now,  $\gamma(d(\oc{t})) = 1$ and assume that   \eqref{eq:lingamma}  is satisfied for some $t$. Note that (\ref{eq:omega}) can be rewritten as:
\begin{equation}\label{eq:proof1}
\begin{split}
&\omega(u_t, d(t-1)) = \\
&D(f(x_t|u_t, x_{t-1}) \| \idf(x_t |u_t, x_{t-1}) )\\
& - \int f(x_t | u_t, x_{t-1})\log(\gamma(d(t)))dx_t.
\end{split}
\end{equation}

Since $f(x_t|u_t, x_{t-1})$ and $\idf(x_t |u_t, x_{t-1})$ are both multivariate normal distributions, and since $x_t$ is $n$-dimensional, it can be shown that:

\begin{equation}\label{eq:proof2}
\begin{split}
D(f(x_t|u_t, x_{t-1}) \| \idf(x_t |u_t, x_{t-1}) ) = & \\
 \frac{1}{2}(\log\frac{|\iSigma_{t,x}|}{|\Sigma_{\Xi}|} - n + \text{tr}(\iSigma_{t,x}^{-1}\Sigma_{\Xi}) \\
 + (\tilde{\mu}_x -\imu_{t,x})^T \iSigma_{t, x} (\tilde{\mu}_x -\imu_{t,x}))& 
 \end{split}
\end{equation}

where $\tilde{\mu}_x = A_t x_{t-1} +B_t u_t$.  Since \eqref{eq:lingamma} is satisfied at time $t$, then computing the second term of the right hand side of (\ref{eq:proof1}) yields 

\begin{equation}\label{eq:proof3}
\begin{split}
\mathbb{E}[\log(\gamma(d(t)))] &= C - \frac{1}{2} \mathbb{E}[(L_t X_t + M_t)^T \Sigma_{t, \gamma}^{-1} (L_t X_t + M_t)] \\
&= C' - \frac{1}{2} (L_t \tilde{\mu}_x + M_t)^T \Sigma_{t,\gamma}^{-1} (L_t \tilde{\mu}_x +M_t)
\end{split}
\end{equation}
where $C' = C - \frac{1}{2} \text{tr}(\Sigma_{\gamma, t}^{-1}L_t\Sigma_{\Xi} L_t^T)$. It can be shown that, by combining (\ref{eq:proof2}) and (\ref{eq:proof3}), by completing the squares and normalizing by $\gamma(d(t-1))$, gives  \eqref{eq:sig-om-t} and \eqref{eq:mu-om-t}. 

Moreover, combining (\ref{eq:proof3}) with \eqref{eq:out} and \eqref{eq:gamma} gives us  \eqref{eq:sig-u-t}, \eqref{eq:mu-u-t}, \eqref{eq:l-t}, \eqref{eq:sig-g-t}, \eqref{eq:m-t}. $\square$\newline

\begin{rem}
In Proposition \ref{lem:full-linear} we considered a model described by $f(x_t| u_t,x_{t-1}) = \N\left(A_tx_{t-1}+B_tu_t,\Sigma_{\Xi}\right)$. It is known that a random vector, $x_t \sim  \N\left(A_tx_{t-1}+B_tu_t,\Sigma_{\Xi}\right)$ is given by the stochastic equation
$$
x_t = A_t x_{t-1} + B_t u_t + \xi_t,
$$ 
where $A_t$ is the time-varying $n \times n$ state matrix, $B_t$ is the time-varying $n \times m$ control matrix and $\xi_t \sim \N(0,\Sigma_{\Xi})$ is some noise. \newline
\end{rem}

\begin{rem}
The recursive equations given in Proposition \ref{lem:full-linear} imply that $x_t$ is distributed following a multivariate normal. In turn, this allows to explicitly compute the  mean vector and the covariance matrix for the pdf describing $x_t | d(t-1)$.\newline
\end{rem}

\subsection*{Example: the optimal control strategy}
\begin{figure}[thbp]
	\includegraphics[width=0.5\textwidth]{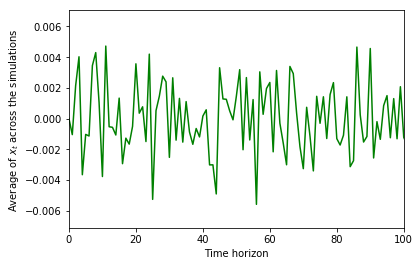}
	\caption{Time evolution of $x_t$ averaged across $1e5$ simulations.}\label{fig:yt}
\end{figure}

\begin{figure}[thbp]
	\includegraphics[width=0.5\textwidth]{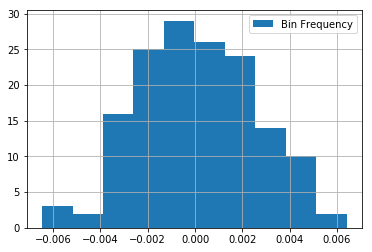}
	\caption{Statistical distribution of $x_t$ obtained from the results of $1e5$ simulations.}\label{fig:bin}
\end{figure}

We now illustrate the effectiveness of the control strategy of Proposition \ref{lem:full-linear} via a representative example. In particular, we consider a single input-single output system. The system model is given by:
\begin{equation}\label{eqn:sys_example}
f(x_t| u_t,x_{t-1}) = \N\left(ax_{t-1}+bu_t,\Sigma_{\Xi}\right),
\end{equation}
where $a = 1.27$, $b = 0.04$, $\Sigma = 0.6$. The ideal distributions are instead:
\begin{equation*}
\begin{split}
\idf(x_t | u_t, x_{t-1}) & = \mathcal{N}(0, \iSigma_{t,x}), \\  
\idf(u_t|d(t-1)) & = \mathcal{N}(0, \iSigma_{t, u})
\end{split}
\end{equation*}
where $\iSigma_{t,x} = 0.2$ and $\iSigma_{t, u} = 0.4$. Finally, we set the time horizon to $\oc{t} =100$. Note that the uncontrolled system is unstable ($a$ is indeed greater than $1$). However, the above result implies that the optimal admissible policy of Proposition \ref{lem:full-linear} is able to stabilize the system and, moreover, ensure that $x_t$ will be distributed in accordance to the specified pdf $\idf(x_t | u_t, x_{t-1})$. This is confirmed in Figure \ref{fig:yt}. The figure has been obtained by running $1e5$ simulations and then by averaging, for each $t$ in the time horizon, the results across all the realizations. Finally, in Figure \ref{fig:bin} we further characterize the statistical distribution of $x_t$ by binning the variable, in order to approximately {\em visualize} the underlying pdf.

\subsection{Link with the Linear Quadratic Regulator}

In the special case when $\imu_{t, x} = \imu_{t, u} = 0$, the recursive equations in Proposition \ref{lem:full-linear} become the well-known update equations as those of the Linear Quadratic Regulator when the following cost function is minimized:
$$
J = x_{\oc{t}}^T\iSigma_{\oc{t}, x}^{-1}x_{\oc{t}} + \sum_{t = 0, \oc{t}-1} x_t^T \iSigma_{t, x}^{-1} x_t + u_t^T \iSigma_{t, u}^{-1} u_t.
$$ 
Therefore all the results for the LQR apply in this set-up and the infinite horizon case will involve solving the same Algebraic Ricatti Equation. Interestingly, this allows to recover the Certainty Equivalence Principle from the Linear Quadratic-Gaussian control \cite{Ste_84}.

Additional considerations can be drawn on the mean and variance of $x_t$. For the sake of notational ease, assume that the control and the state matrices are constant over time and note that the time evolution of $x_t$ can be seen as the solution of the stochastic equation
\begin{equation}\label{eq:LTI}
x_t = Ax_{t-1} + B\mu_{t,u} + \xi_t + \xi_{u,t},
\end{equation}
where $\xi_{u,t} \sim \mathcal{N}(0, \Sigma_{t, u})$. Then:
\begin{subequations}
\begin{equation}\label{eq:mu-x-t}
\begin{split}
\mu_{t, x} & = (A - B\Sigma_{t,u} B^T \Sigma_{t,\omega}^{-1} A)\mu_{t-1, x} \\
& = K\mu_{t-1, x},
\end{split}
\end{equation}    
\begin{equation}\label{eq:-x-t}
\Sigma_{t, x} = \Sigma_{\Xi} + K\Sigma_{t-1, x} K^{T} + B \Sigma_{t, u} B^T.
\end{equation}
\end{subequations}

In particular, since the controller is randomized, from (\ref{eq:-x-t}) we see that the variance of the controller (which is still the optimal controller in the KL-Divergence metric) affects the variance of the state variable. This is the major difference between the randomized control and the LQG or LQR control algorithms. In such algorithms, the control variable is chosen deterministically once $x_{t-1}$ is known. This observation leads us to the following proposition.\newline

\begin{Proposition}
	\label{lem:stab}
Assume that the probabilistic control given by the equations in Proposition \ref{lem:full-linear} is used to control (\ref{eq:LTI}). Then, the conditions for stability of the closed loop system (both in variance and in mean) are the same conditions as in the ones of the LQR. Namely, that all the eigenvalues of the matrix $S = (A-B\Sigma_{t,u}B^T\Sigma_{t,\omega}^{-1}A)$ lie inside the unit circle.\newline
\end{Proposition}
The proof of this result is omitted here for the sake of brevity.

\subsection{Analysis of model error propagation}\label{sec:model-error}

Now, we turn our attention to study how sensitive the control algorithm in Proposition \ref{lem:full-linear} is with respect to errors in the models. In doing so, we consider the infinite horizon setting, i.e. when $\oc{t}\rightarrow+\infty$. Also, we assume that $\iSigma_{t,u}$, $\iSigma_{t,x}$ are constant and $\imu_{t,u} = \imu_{t,x} = 0$. Moreover, in the context of this analysis, the {\em real} statistical model for the system is given by:
$$
f(x_t | u_t,x_{t-1}) = \N(Ax_{t-1}+Bu_t,\Sigma_{\xi}),
$$
while the control algorithm of Proposition \ref{lem:full-linear} does not have access to the real model, but only to its approximation. Specifically, the control algorithm has access to the {\em approximated} system model
\begin{equation}\label{eq:approximate_model}
\tilde f(x_t | u_t,x_{t-1}) = \N(\tilde Ax_{t-1}+\tilde Bu_t,\tilde\Sigma_{t,\xi}).
\end{equation}
\begin{rem}
Situations where the control algorithm has only access to a statistical model that can only approximate the real model of the system naturally arise in many applications. For example, it is known that, in the context of deep learning, neural networks \cite{Goodfellow-et-al-2016} can be used to approximate the pdf of the underlying process being classified. In other applications involving iterative learning, the presence of bias in the training datasets can lead to biased models which could considerably affect the outcomes of the learning process, see e.g. \cite{Hasselt:2016:DRL:3016100.3016191}. With the analysis present here we seek to investigate the impact of model approximations on the control algorithm of Proposition \ref{lem:full-linear}.\newline
\end{rem}

Given the set-up of this Section, the recursive equations for the control algorithm of Proposition \ref{lem:full-linear} are considerably simplified. This is formalized with the following proposition.\newline

\begin{Proposition}\label{lem:infinite_horizon}
Consider the control algorithm of Proposition \ref{lem:full-linear} and assume that: (i) $\oc{t}\rightarrow+\infty$; (ii) $\iSigma_{t,u}$, $\iSigma_{t,x}$ are constant; (iii) $\imu_{t,u} = \imu_{t,x} = 0$; (iv) the approximate model used by the algorithm is given in (\ref{eq:approximate_model}). Then, the randomized control strategy is given by the  system of recursive equations
\begin{subequations}
\begin{equation}\label{eq:s-u-t-sim}
\Sigma_{u}^{-1} = \iSigma_{u}^{-1} + \tilde{B}^T \Sigma_{\omega}^{-1} \tilde{B},
\end{equation}    
\begin{equation}\label{eq:mu-u-t-sim} 
\mu_{t,u} = -\Sigma_{u}\tilde{B}^T \Sigma_{\omega}^{-1}\tilde{A} x_{t-1},
\end{equation}
\begin{equation}\label{eq:s-g-t-sim} 
\Sigma_{\gamma}^{-1} = \Sigma_{\omega}^{-1} (\Sigma_{\omega} - \tilde{B}\Sigma_{u}\tilde{B}^T)\Sigma_{\omega}^{-1},
\end{equation}
\end{subequations}
where $\Sigma_{\omega}^{-1}$ solves the algebraic Riccati equation
\begin{equation}\label{eq:s-t-o-sim}
\begin{split}
	\Sigma_{\omega}^{-1} &= \iSigma_{x}^{-1} + \tilde{A}^T \Sigma_{\omega}^{-1} \tilde{A} \\
	& - \tilde{A}^T \Sigma_{\omega}^{-1} \tilde{B}(\iSigma_{u}^{-1} + \tilde{B}^T \Sigma_{\omega}^{-1}\tilde{B})^{-1}\tilde{B}^{T}\Sigma_{\omega}^{-1}\tilde{A}.
\end{split}
\end{equation}\newline
\end{Proposition}

The proof of the above result can be obtained via direct inspection, by using the assumptions of Proposition \ref{lem:infinite_horizon} to simplify the equations given in Proposition \ref{lem:full-linear}.\newline

From the above proposition note that the randomized policy can be computed by solving a Riccati equation. That is, in order to study  stability of the system we need to study the eigenvalues of the matrix $(A - B\Sigma_{u}\tilde{B}^T\Sigma_{\omega}\tilde{A})$. Now, consider the case where the system controlled by the algorithm of Proposition \ref{lem:infinite_horizon} is stable. That is, we assume that $\rho(A-B\Sigma_uB^T\Sigma_{\omega}A)<1$ where $\rho(M)$ is the spectral radius of the matrix $M$. Clearly, when  $A \ne \tilde{A}$ and $B \ne \tilde{B}$ then such a condition can be violated. We are interested in identifying for which set of parameters, $\tilde A$ and $\tilde B$, stability is preserved. This leads to the following:\newline

\begin{Definition}
The region of convergence of a system controlled by the algorithm of Proposition \ref{lem:infinite_horizon}  is the region of the parameters space 
$$
K_{A,B} = \{(\tilde{A},\tilde{B}): \rho( A - B\Sigma_u\tilde{B}^T\Sigma_{\omega}\tilde{A}) < 1\}.
$$\newline
\end{Definition}

Since $\rho(\cdot)$ is continuous, $K_{A, B}$ is an open set that includes at least one point in the space parameter, ($A, B$). Therefore it must contain a neighborhood of that point. Furthermore, the region of convergence must necessarily be bounded in the $\tilde{A}$ direction for a fixed $\tilde{B} \neq 0$ (if $\rho(A) < 1$ then $(\tilde{A}, 0) \in K_{A, B}$ for all $\tilde{A}$). The behavior with respect to $\tilde{B}$ is more difficult to analyze in the matrix case due to the implicit definition of $\Sigma_{\omega}^{-1}$. We can, however, find the boundary of $K_{A,B}$ using an algorithm for continuation of an implicitly defined manifold, for example \cite{henderson_multiple_2002}. These algorithms need an initialization point on the boundary. In this case, a point where $ \rho(A - B\Sigma_u\tilde{B}^T\sigma_{\omega}\tilde{A}) = 1$ or close enough. The point can be found using any root finding algorithm. 

The above observations imply that the region of convergence can be numerically computed  via continuation methods if some points of the parameter space are known to belong to the region.  Suppose we have some theoretical knowledge on $A$ and $B$, namely, we are given a subset $\mathcal{K}$ of the parameter space where we know the system is stable. We are interested now in computing the set:
$$
K_{\mathcal{K}} = \cap_{(A, B)\in\mathcal{K}} K_{A,B}.
$$ 
Alternatively, one can define this set as
$$
K_{\mathcal{K}} = \{(\tilde{A},\tilde{B}): \max_{(A,B)\in\mathcal{K}} \rho(A - B\Sigma_u\tilde{B}^T\Sigma_{\omega}\tilde{A}) < 1\}.
$$ 
Again, the function $\max_{(A,B)\in\mathcal{K}} \rho(A - B\Sigma_u\tilde{B}^T\sigma_{\omega}\tilde{A})$ is continuous and therefore we can compute the boundary $\max_{(A,B)\in\mathcal{K}} \rho(A - B\Sigma_u\tilde{B}^T\sigma_{\omega}\tilde{A}) = 1$ using the algorithm in \cite{henderson_multiple_2002}. Thus, we have the following proposition, stating a condition ensuring {\em computability} of the convergence region.\newline

\begin{Proposition}\label{lem:stab_region}	
Given theoretical knowledge of a compact set $\mathcal{K}$, we can compute the set $K_{\mathcal{K}}$, guaranteeing that any $\tilde{A},\tilde{B} \in K_{\mathcal{K}}$ will make the system stable.
\end{Proposition}

\subsection*{Example (continued): effects of model errors and convergence region}
\begin{figure}[thbp]
	\includegraphics[width=0.5\textwidth]{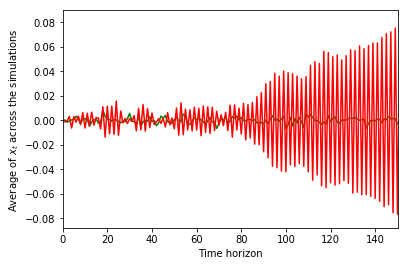}
	\caption{Time evolution of $x_t$ averaged across $1e5$ simulations when the approximate model (\ref{eqn:appr_model_example}) is used by the algorithm of Proposition \ref{lem:infinite_horizon}. In green: $x_t$ when the real model is used by the control algorithm. In red: $x_t$ when the approximate model is instead  used. Colors online}\label{fig:yt_error}
\end{figure}

\begin{figure}[thbp]
	\includegraphics[width=0.5\textwidth]{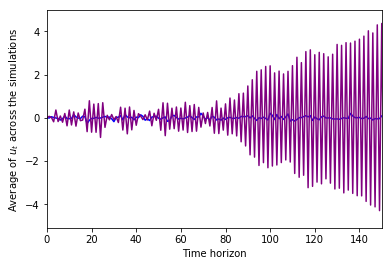}
	\caption{Time evolution of $u_t$  averaged across $1e5$ simulations when the approximate model (\ref{eqn:appr_model_example}) is used by the algorithm of Proposition \ref{lem:infinite_horizon}. In blue: $u_t$ when the real model is used by the control algorithm. In purple: $u_t$ when the approximate model is instead  used. Colors online}\label{fig:ut_error}
\end{figure}

Consider, again, the system model (\ref{eqn:sys_example}). This time, we consider the infinite horizon control problem and assume that the control algorithm does not have a complete knowledge of the model but rather it only has access to an approximate version, namely:
\begin{equation}\label{eqn:appr_model_example}
\tilde f(x_t | u_t,x_{t-1}) = \N(ax_{t-1}+\tilde bu_t,\Sigma_{\xi}),
\end{equation}
where $\tilde b = 0.02$. That is, the only difference of the approximate model with respect to the real system model lies in the parameter $b$. As shown in Figure \ref{fig:yt_error} and Figure \ref{fig:ut_error} this mismatch leads to significant changes in the behavior of the closed loop system, when the control algorithm of Proposition \ref{lem:infinite_horizon} is used.

\begin{figure}[h]
\begin{tabular}{cc}
\includegraphics[width=0.21\textwidth]{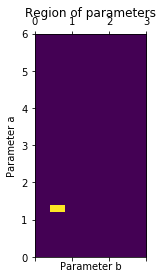} & \includegraphics[width=0.21\textwidth]{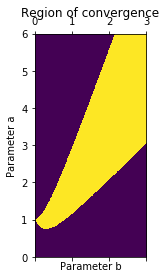}
\end{tabular}
\caption{For a given region in which the parameters lie, we find a region in which the system converges (in yellow, in both panels). Left panel: initial point of convergence, obtained by solving the Riccati equation. Right panel: region of convergence estimated via continuation methods \cite{henderson_multiple_2002}.}\label{fig:region_convergence}
\end{figure}

Since the system is scalar, we were able to solve explicitly the Riccati equation in order to estimate the region of convergence of the control algorithm. In particular, once we found one point in the parameter space for which the closed loop system converged (left panel of Figure \ref{fig:region_convergence}) we were able to estimate the region $K_{\mathcal{K}}$, in accordance to Proposition \ref{lem:stab_region}.

\subsection{Safety of the closed loop system}\label{sec:safety}

Given the set-up of Section  \ref{sec:model-error}, we now consider {\em safey} of the closed loop system, when the control algorithm of Proposition \ref{lem:infinite_horizon} is used to control a system of interest, for which only the approximate model (\ref{eq:approximate_model}) is known. Intuitively, inspired by e.g. \cite{safety_18}, we define the safety constraints of the closed-loop system in terms of a {\em safety region}, i.e. a region of the state space that is, with {\em high probability}, forward invariant. Essentially, if the initial conditions of the system start in the safety region, then this region, with high probability, is never left by the system trajectories. In this paper we are interested in characterizing the safety condition as a function of $\tilde A$ and $\tilde B$. We consider the case where the initial conditions, $x_0$, of the system are safe and specify the safety region in terms of $x_0$. This is formalized with the next definition, where $\| \cdot \|$ denotes Euclidean norm:\newline

\begin{Definition}\label{def:safety}
Let $M >0$, $\delta > 0$ and let $x_0$ be the initial system conditions. The $(M,\delta)$-safety region for the closed loop system is the region of the parameter space
$$
K_{A, B}(M, \delta) = \{\tilde{A},\tilde{B}: \mathbb{P}(\|x_t\| > M\|x_0\|) < \delta \; \forall t\}.
$$\newline
\end{Definition}

\begin{rem}
A major difference of Definition \ref{def:safety} with respect to the one given in \cite{NIPS2017_6692,safety_18} is that the $(M,\delta)$-safety region is a region of the system parameter space, rather than the system state space. \newline
\end{rem}

The following proposition is an analogous of Proposition \ref{lem:stab_region} and ensures the computability of the safety region.\newline

\begin{Proposition}
Assume that the initial conditions, $x_0$, of the system are safe. Then, the safety region given in Definition \ref{def:safety} can be computed for any $\delta >0$ and $M > 0$.\newline
\end{Proposition}

Essentially, since the distribution of $x_t$ is controlled, the set $K_{A, B}(M, \delta)$ is, again, computable via the algorithmic continuation approach of e.g. \cite{henderson_multiple_2002}. \newline

\begin{rem}
In principle, computing the sets $K_{\mathcal{K}}$ and $K_{\mathcal{K}}(M, \delta)$ allows to select appropriate $\tilde{B}, \tilde{A}$ under uncertainty in the model $A, B$. This might be of use in a system identification scenario, when the parameters of the identified system are uncertain but bounds on uncertainty are known.
\end{rem}

\subsection{Integrating iterative learning of the model}

Suppose now that we are in a scenario in which the system model is learned via an iterative learning process (see e.g. \cite{8039204} for an application in the context of model predictive control). From the viewpoint of the probabilistic control analyzed in this paper, we are not interested in designing the learning strategy. Rather, we model the learning process by having $\tilde{A}_t, \tilde{B}_t$ such that:
$$
\|A-\tilde{A}_t\| \leq \epsilon_1 + \epsilon_2^t,
$$ 
and
$$
\|B-\tilde{B}_t\| \leq \epsilon_3 + \epsilon_4^t,
$$ 
where $\varepsilon_2$ and $\varepsilon_4$ are  smaller than one. Intuitively, with the above equations we model the case where the learning process allows to update the matrices $\tilde{A}_t$ and $\tilde{B}_t$ so that such matrices get closer to the real ones. Also, note that we consider the non-ideal case, where we allow for the learning process to convergence to matrices that are close to the real ones. This is modeled by the fact that, when $t\rightarrow+\infty$, we have $\|A-\tilde{A}_t\| \leq \epsilon_1$ and $\|B-\tilde{B}_t\| \leq \epsilon_3$. 

Given the scenario illustrated above, it can be shown (the proof will be presented elsewhere) that the (stable) control policy can be computed as follows:
\begin{enumerate}
\item assume that, at each time step, $t$, the matrices $\tilde{A}_t$ and $\tilde{B}_t$ are equal to the real matrices;
\item assume that such matrices will not change over time;
\item apply \eqref{eq:s-u-t-sim} and \eqref{eq:mu-u-t-sim} on such matrices.
\end{enumerate}

Finally, the introduction of a learning mechanism modeled as described above, also improves the characterization of the safety region. In particular, it can be shown that the region can be computed and, moreover, it is possible to compute the set of initial conditions for which the system is safe. The numerical methods, which are still based on the use of continuation algorithms,  are omitted here and will be presented elsewhere.

\subsection*{Example (continued): safety region}
\begin{figure}
 \includegraphics[width=0.5\textwidth]{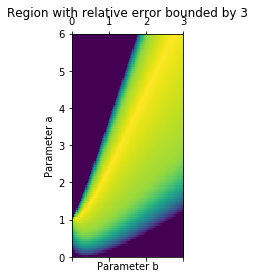}
 \caption{Safety region with $M= 3$ and $\delta =0.1$. The color scale indicates the maximum of $\| x_t \|$: the region in purple does not belong to the safety region. Colors online.}\label{relative_error_1}
\end{figure}

Finally, we considered again the approximate model in (\ref{eqn:appr_model_example}). This time, we assume that the parameter $\tilde b$ changes over time and it indeed converges exactly to $b$. In this case, we were able to compute the $(3,0.1)$-safety region of the system, as shown in Figure \ref{relative_error_1}.

\section{Conclusions}
We presented a sensitivity and a safety analysis for the so-called fully probabilistic control scheme. In particular, after presenting the sensitivity analysis, we focused on characterizing the convergence region of the closed loop system and introduced a safety analysis. Finally, we discussed how the introduction of  learning mechanisms can be beneficial for both convergence and safety.

\section*{Acknowledgments}
BGP was partially supported by the NSF grant number 1514606. GR was supported in part by a research grant from Science Foundation Ireland (SFI) under Grant Number 16/RC/3872.

\bibliographystyle{IEEEtran}
\bibliography{fpd_robust}           

\end{document}